\documentclass[12pt]{amsart}
\usepackage{amssymb}
\usepackage{amsthm}
\usepackage{euscript}

\theoremstyle{plain}

\newtheorem{thm}[subsection]{Theorem}
\newtheorem{prop}[subsection]{Proposition}

\begin{document}
\title[Burnside's theorem in the setting of general fields]{Burnside's theorem in the setting of general fields}
\author{Heydar Radjavi and Bamdad R. Yahaghi}

\address{Department of Pure Mathematics, University of Waterloo, Waterloo, Ontario, Canada N2L 3G1 \newline
\indent Department of Mathematics, Faculty of Sciences, Golestan University, Gorgan 19395-5746, Iran}
\email{hradjavi@uwaterloo.ca, bamdad5@hotmail.com,\newline  bamdad@bamdadyahaghi.com}


\bibliographystyle{plain}

\begin{abstract}
We extend a well-known theorem of Burnside in the setting of general fields as follows:  for a general field $F$ the matrix algebra $M_n(F)$ is the only algebra in $M_n(F)$  which is spanned by an irreducible semigroup of triangularizable matrices. In other words, for a semigroup of triangularizable matrices with entries from a general field irreducibility is equivalent to absolute irreducibility. As a consequence of our result we prove a stronger version of a theorem of Janez Bernik.
\end{abstract}

\maketitle

\bigskip

\begin{section}
{\bf Introduction}
\end{section}

\bigskip

A version of a celebrated theorem of Burnside \cite[Theorem on p. 433]{Bu} asserts that for an algebraically closed field $F$, the matrix algebra $M_n(F)$ is the only  algebra in $M_n(F)$  which is spanned by an irreducible semigroup of matrices. We prove a counterpart of Burnside's Theorem  in the setting of general fields as follows:  for a general field $F$, the matrix algebra $M_n(F)$ is the only  algebra in $M_n(F)$  which is spanned by an irreducible semigroup of triangularizable matrices.  In other words, for a semigroup of triangularizable matrices with entries from a general field irreducibility  is equivalent to absolute irreducibility.

Throughout, $F$ and $K$ stand for fields and $F$ is a subfield of $K$. We view the elements of the matrix algebra $M_n(F)$ as linear transformations acting on the left of $F^n$, the vector space of all $n \times 1$ column vectors with entries from $F$. A family $\mathcal{F}$  in $M_n(F)$ is said to be irreducible if the orbit of any nonzero $ x \in F^n$ under the algebra generated by $\mathcal{F}$, denoted by ${\rm Alg}(\mathcal{F})$, is $  F^n$. When $ n > 1$, this is easily seen to be equivalent to the lack of nontrivial invariant subspaces for the family $\mathcal{F}$; the trivial spaces being $\{0\}$ and $  F^n$. Reducible, by definition, means not irreducible. Absolute irreducibility means irreducibility over any field extension of the ground field $F$. It follows from Burnside's Theorem that a family $\mathcal{F}$  in $M_n(F)$ is absolutely irreducible if and only if  ${\rm Alg}(\mathcal{F})= M_n(F)$. On the opposite side of irreducibility, we have the notion of triangularizability. More precisely, a family  $\mathcal{F}$  in $M_n(F)$ is called triangularizable if there exists a maximal  chain
$$\{0\} = \mathcal{M}_0 \subset  \mathcal{M}_1 \subset  \cdots \subset \mathcal{M}_n= F^n$$
 of subspaces of $F^n$
with $  \mathcal{M}_i$'s being invariant under the family  $\mathcal{F}$. Any such chain is called a triangularizing chain of subspaces for the family  $\mathcal{F}$. It is a standard observation that a family  $\mathcal{F}$  in $M_n(F)$ is  triangularizable if and only if there exists a basis for $  F^n$, called a triangularizing basis, relative to which every element of the family has an upper triangular matrix. This occurs if and only if there exists an invertible matrix $P \in M_n(F)$ such that the family $ P^{-1} \mathcal{F} P$  consists of upper triangular matrices.

\bigskip

\begin{section}
{\bf Main Results}
\end{section}

\bigskip

The following can be thought of as an extension of Burnside's Theorem to general fields. In fact, it extends  \cite[Theorem 2.3]{Y3} to arbitrary fields.

\bigskip

\begin{thm} \label{2.1}
Let $ n \in \mathbb{N}$,  $F$ be a field, and $ \mathcal S$ a semigroup of triangularizable matrices in $ M_n(F)$. Then the semigroup  $ \mathcal S$ is irreducible iff it is absolutely irreducible.
\end{thm}

\bigskip

\noindent {\bf Proof.} The ``if" implication is trivial. We prove ``the only if" implication. Let $ \mathcal S$ be a semigroup of triangularizable matrices in $ M_n(F)$ and let   $\mathcal{A}= {\rm Alg} (\mathcal{S})$ denote the algebra generated by $\mathcal{S}$ and $r$ the minimal nonzero rank present in  $\mathcal{A}$. As shown in \cite{Y1} (the remark following Theorem 2.9 of that paper),  $r$ divides $n$   and    $\mathcal{A}$ is simultaneously similar to $ M_{n/r} ( \Delta)$, where  $ \Delta$ is an irreducible division algebra in $M_r(F)$, which is necessarily  of dimension $r$. This in particular implies that the algebra $\mathcal{A}$ is simple (and semisimple). We prove the assertion by showing that $ r = 1$. Let $F_c$ denote the algebraic closure of $F$ and view  $\mathcal{A}$ as a simple $F$-algebra in $M_n(F_c)$. Apply a simultaneous similarity to put $\mathcal{A} \subseteq M_n(F_c)$ in block upper triangular form so that the number $k$ of the diagonal blocks is maximal and hence each diagonal block is absolutely irreducible. Note that the diagonal blocks are all nonzero because $\mathcal{A}$ contains the identity matrix.  If necessary,  using \cite[Theorem 1.1]{B} and applying a simultaneous similarity, we may assume that each diagonal block of $\mathcal{A}$  is the full matrix algebra $M_{n_i}(F)$ for some $n_i \in \mathbb{N}$ ($ 1 \leq i \leq k)$. For each $ 1 \leq i \leq k$, let $\mathcal{A}_i = M_{n_i}(F) $ denote the $i$-th diagonal block of $\mathcal{A}$.  In view of the simplicity of $\mathcal{A}$, the mapping $ \phi_{ij} : \mathcal{A}_i  \longrightarrow \mathcal{A}_j$ defined by $ \phi_{ij}(A_i) = A_j$ is a well-defined nonzero homomorphism of $F$-algebras whose inverse $ \phi_{ji} : \mathcal{A}_j  \longrightarrow \mathcal{A}_i$  is also a homomorphism of $F$-algebras. Thus $n_i = n_j= n/k$, and hence $\mathcal{A}_i  = \mathcal{A}_j =M_{n/k} (F)$ for each $ 1 \leq i , j \leq k$. Since $F$-algebra automorphisms of $M_{n/k}(F)$ are  all inner, again if necessary applying another simultaneous similarity, we may assume that the diagonal blocks of each element of $A \in \mathcal{A}$ are all of the size $n/k$ and equal. It thus follows from the Noether-Skolem Theorem, \cite[p. 39]{D}, that the $F$-algebra $ \mathcal{A}$ is similar to the $k$-fold inflation of the $F$-algebra $ M_{n/k}(F)$ because it is isomorphic to it. This in particular implies $ r= k$ and  $ \dim \mathcal{A} = \dim  M_{n/r} ( \Delta)= \dim M_{n/k}(F)$. This clearly yields $ k^2 = k$, and hence $ r = k = 1$, proving the assertion.
\hfill \qed

\bigskip

We used the main result of Bernik  \cite[Theorem 1.1]{B} to prove Theorem \ref{2.1}. The following shows that Theorem \ref{2.1} implies a stronger version of the main result of  \cite{B}. Therefore, any independent proof of Theorem \ref{2.1} would provide a new proof of Bernik's result. (Thus it should be observed that the following proof does not use Bernik's Theorem.)

\bigskip

\begin{thm} \label{2.2}
Let $ n \in \mathbb{N}$, $F$ and $K$ be fields with  $F \leq K$, and $ \mathcal S$ an irreducible semigroup of triangularizable matrices in $ M_n(K)$ with spectra in $F$. Then ${\rm Alg}_F (\mathcal S)$ is similar to $ M_n(F)$ over $M_n(K)$.
\end{thm}

\bigskip

\noindent {\bf Proof.} Since irreducibility implies absolute irreducibility for semigroups of triangularizable matrices, we see that $ \{ 0\} \not= {\rm tr}(\mathcal{S}) \subseteq F$. It thus follows from \cite[Corollary 2.8]{Y1} and absolute irreducibility of  $ \mathcal S$ that $\dim_F {\rm Alg}_F (\mathcal S) = \dim_K {\rm Alg}_K (\mathcal S) = n^2$. Let $ \mathcal{B} \subseteq \mathcal{S}$ be a basis for $\mathcal{A}:=  {\rm Alg}_F (\mathcal S)$ and for $A \in \mathcal A$, $L_A : \mathcal{A} \longrightarrow  \mathcal{A}$, defined by $L_A(B) = AB$, be the linear operator of left multiplication by $A$. It is plain that the mapping $ \phi: \mathcal{A} \longrightarrow M_{n^2}(F)$ defined by $\phi(A) = [L_A]_\mathcal{B}$, where $[L_A]_\mathcal{B}$ denotes the matrix representation of $L_A$ with respect to the basis $\mathcal{B}$, is an embedding of the $F$-algebra $\mathcal{A}$ in $M_{n^2}(F)$. Clearly, $ \phi(\mathcal{A}) $ is a simple subalgebra of $M_{n^2}(F)$. Apply a simultaneous similarity to put $\phi(\mathcal{A}) \subseteq M_{n^2}(F)$ in block upper triangular form so that the number $k$ of the diagonal blocks is maximal and hence each diagonal block is irreducible. Note that the diagonal blocks are all nonzero because $\mathcal{A}$ contains the identity matrix.  Also note that $\mathcal{A}=  {\rm Alg}_F (\mathcal S)$ and $ \mathcal{S}$ consists of triangularizable matrices.  But irreducibility implies absolute irreducibility for semigroups of triangularizable matrices. Thus, each diagonal block of $ \phi(\mathcal{A}) $   is an absolutely irreducible $F$-algebra in $M_{n_i}(F)$,  and hence is equal to the full matrix algebra $M_{n_i}(F)$. From this point on, an argument almost identical to that of the proof of Theorem \ref{2.1} shows that the $F$-algebra $ \phi(\mathcal{A})$ is similar to the $k$-fold inflation of the $F$-algebra $ M_{n^2/k}(F)$ for some $k$ dividing $n^2$. This in particular gives
$$ n^2 = \dim_F \mathcal{A} = \dim_F \phi(\mathcal{A}) = \dim_F  M_{n^2/k}(F) = n^4/k^2, $$
which in turn implies $ n = k$. Consequently, the $F$-algebra  $\mathcal{A}$ is isomorphic to the $1$st block diagonal of  $ \phi(\mathcal{A})$, which is $ M_n(F)$. It thus follows from the Noether-Skolem  Theorem, \cite[p. 39]{D},  that $\mathcal{A}$ is similar to  $ M_n(F)$, which is the desired result.
\hfill \qed

\bigskip

\noindent  {\bf Remark.} With this theorem at our disposal, we can prove the counterparts of \cite[Theorem B on p. 99]{K} and \cite[Theorem 1]{RR2} for semigroups of triangularizable matrices in $ M_n(K)$ with spectra in $F$, see \cite[Theorems 2.7 and 2.8]{Y3}

\bigskip

An extension of Burnside's Theorem was proved in \cite[Theorems 2.1-2]{RY} as follows: for an $n > 1$ and a  finite field, or more generally a quasi-algebraically closed field $F$, $M_n(F)$ is the only irreducible algebra in $M_n(F)$ that, as a vector space over $F$, is spanned by triangularizable matrices in $M_n(F)$.  The following proposition answers a question left open in  \cite[Remark 1 following  Theorem 2.1]{RY}.  More precisely, it shows that the theorem does not hold for general fields, e.g., for the real field because $ M_n(\mathbb{H})$, viewed as a proper irreducible subalgebra of $M_{4n}(\mathbb{R})$, is spanned by the identity matrix and nilpotents as a vector space over $\mathbb{R}$.

\bigskip

\begin{prop} \label{2.3}
Let $ n \in \mathbb{N}$ with $ n > 1$ and  $ \mathbb{H}$ denote the division ring of quaternions. Then $ M_n(\mathbb{H})$ is spanned by $I$, the identity matrix, and nilpotents as a vector space over $\mathbb{R}$.
\end{prop}

\bigskip

\noindent {\bf Proof.} It suffices to prove the assertion for $M_2( \mathbb{H})$. Since
$\left(\begin{array}{cc}
0 & p \\
0& 0
\end{array}\right)$,
$\left(\begin{array}{cc}
0 & 0 \\
q& 0
\end{array}\right)$,
and
$\left(\begin{array}{cc}
p & p \\
-p& -p
\end{array}\right)$
are all nilpotents for all  $ p, q \in \mathbb{H}$, we only have to show that
$\left(\begin{array}{cc}
p & 0 \\
0& p
\end{array}\right)$
is spanned by $I$ and nilpotents for all $ p \in \mathbb{H}$. Let $ p = a + b i + c j + d k$, where $ a, b, c, d \in \mathbb{R}$. Thus it suffices to show that
$\left(\begin{array}{cc}
i & 0 \\
0& i
\end{array}\right)$,
$\left(\begin{array}{cc}
j & 0 \\
0& j
\end{array}\right)$,  and
$\left(\begin{array}{cc}
k & 0 \\
0& k
\end{array}\right)$
are all in the desired span. Now
$$ \left(\begin{array}{cc}
i & 0 \\
0& i
\end{array}\right)=
\left(\begin{array}{cc}
i & j \\
-j& i
\end{array}\right)+
\left(\begin{array}{cc}
0 & -j \\
0& 0
\end{array}\right)+
\left(\begin{array}{cc}
0 & 0 \\
j& 0
\end{array}\right).$$
But
$ \left(\begin{array}{cc}
i & j \\
-j& i
\end{array}\right)^2 = 0$. This completes the proof.
\hfill \qed


\end{document}